\documentclass[preprint,11pt]{elsarticle}
\usepackage[latin1]{inputenc}
\usepackage[english]{babel}
\usepackage{geometry}
\usepackage{amsmath}
\usepackage{amssymb}
\usepackage{oldgerm}
\usepackage{amscd}
\usepackage{hyperref}

\let\today\relax

\newtheorem{Def}{Definition}[section]
\newtheorem{Prop}{Proposition}[section]
\newtheorem{Theorem}{Theorem}[section]
\newtheorem{Lemma}{Lemma}[section]
\newtheorem{Cor}{Corollary}[section]
\newcommand{\Proof}{\noindent{\bfseries Proof : }} 
\newdefinition{Oss}{Remark}
\newcommand{\Q}{\mathbb{Q}}
\newcommand{\N}{\mathbb{N}}
\newcommand{\Z}{\mathbb{Z}}
\newcommand{\IZ}{\textnormal{Int($\Z$)}}

\makeatletter
\def\ps@pprintTitle{%
     \let\@oddhead\@empty
     \let\@evenhead\@empty
     \def\@oddfoot{\reset@font\hfil\thepage\hfil%
     \llap{\footnotesize\itshape\today}}
     \let\@evenfoot\@oddfoot}
\makeatother

\begin{document}

\begin{frontmatter}

\title{Primary decomposition of the ideal of polynomials whose fixed divisor is divisible by a prime power\\ \vskip0.3cm
\small{J. Algebra 398 (2014), 227-242, \href{http://dx.doi.org/10.1016/j.jalgebra.2013.09.016}{http://dx.doi.org/10.1016/j.jalgebra.2013.09.016}}}


\author{Giulio Peruginelli\fnref{fn1}}
\ead{peruginelli@math.tugraz.at}
\fntext[fn1]{Old affiliation, from June 2014: Department of Mathematics, University of Padova, Italy.}
\address{Institut f\"ur Analysis und Comput. Number Theory, Technische Universit\"at, Steyrergasse 30, A-8010 Graz, Austria.}


\begin{abstract}
\noindent We characterize the fixed divisor of a polynomial $f(X)$ in $\Z[X]$ by looking at the contraction of the powers of the maximal ideals of the overring $\IZ$ containing $f(X)$. Given a prime $p$ and a positive integer $n$, we also obtain a complete description of the ideal of polynomials in $\Z[X]$ whose fixed divisor is divisible by $p^n$ in terms of its primary components.
\end{abstract}


\end{frontmatter}

\vskip0.2cm
\begin{flushright}
\small{\emph{to Sergio Paolini, whose teachings and memory I deeply preserve.}}
\end{flushright}

\vskip0.4cm
\section{Introduction}
In this work we investigate the image set of integer-valued polynomials over $\Z$. The set of these polynomials is a ring usually denoted by:
$$\IZ\doteqdot\{f\in\Q[X]\;|\;f(\Z)\subset\Z\}.$$

Since an integer-valued polynomial $f(X)$ maps the integers in a subset of the integers, it is natural to consider the subset of the integers formed by the values of $f(X)$ over the integers and the ideal generated by this subset. This ideal is usually called the fixed divisor of $f(X)$. Here is the classical definition.

\begin{Def}
 Let $f\in\IZ$. The \textbf{fixed divisor} of $f(X)$ is the ideal of $\Z$ generated by the values of $f(n)$, as $n$ ranges in $\Z$:
 $$d(f)=d(f,\Z)=(f(n)|n\in\Z).$$ 
 We say that a polynomial $f\in\IZ$ is \textbf{image primitive} if $d(f)=\Z$.
\end{Def}

It is well-known that for every integer $n\geq1$ we have
\begin{equation*}
d(X(X-1)\ldots(X-(n-1)))=n!
\end{equation*} 
so that the so-called binomial polynomials $B_n(X)\doteqdot X(X-1)\ldots(X-(n-1))/n!$ are integer-valued (indeed, they form a free basis of $\IZ$ as a $\Z$-module; see \cite{CaCh0}).

Notice that, given two integer-valued polynomials $f$ and $g$, we have $d(fg)\subset d(f)d(g)$ and we may not have an equality. For instance, consider $f(X)=X$ and $g(X)=X-1$; then we have $d(f)=d(g)=\Z$ and $d(fg)=2\Z$. If $f\in\IZ$ and $n\in\Z$, then directly from the definition we have $d(nf)=nd(f)$. If cont($F$) denotes the content of a polynomial $F\in\Z[X]$, that is, the greatest common divisor of the coefficients of $F$, we have $F(X)=$cont$(F)G(X)$, where $G\in\Z[X]$ is a primitive polynomial (that is, cont($G$)=1). We have the relation:
$$d(F)=\mbox{cont}(F)d(G).$$
In particular, the fixed divisor is contained in the ideal generated by the content. Hence, given a polynomial with integer coefficients, we can assume it to be primitive. In the same way, if we have an integer-valued polynomial $f(X)=F(X)/N$, with $f\in\Z[X]$ and $N\in\N$, we can assume that $($cont$(F),N)=1$ and $F(X)$ to be primitive.

The next lemma gives a well-known characterization of a generator of the above ideal  (see \cite[Lemma 2.7]{ACaChS}). 
\begin{Lemma}\label{fixdiv}
 Let $f\in\IZ$ be of degree $d$ and set 
 \begin{itemize}
  \item[1)] $d_1=\sup\{n\in\Z\;|\;\frac{f(X)}{n}\in\IZ\}$
  \item[2)] $d_2=GCD\{f(n)\;|\;n\in\Z\}$
  \item[3)] $d_3=GCD\{f(0),\ldots,f(d)\}$ 
 \end{itemize} 
 then $d_1=d_2=d_3$.
\end{Lemma}
Let $f\in\IZ$. We remark that the value $d_1$ of Lemma \ref{fixdiv} is plainly equal to:
$$d_1=\sup\{n\in\Z\;|f\in n\IZ\}.$$
Moreover, given an integer $n$, we have this equivalence that we will use throughout the paper, a sort of ideal-theoretic characterization of the arithmetical property that all the values attained by $f(X)$ are divisible by $n$:
$$f(\Z)\subset n\Z\Longleftrightarrow f\in n\IZ$$
($n\IZ$ is the principal ideal of $\IZ$ generated by $n$).
From 1) of Lemma \ref{fixdiv} we see immediately that if $f(X)=F(X)/N$ is an integer-valued polynomial, where $F\in\Z[X]$ and $N\in\N$ coprime with the content of $F(X)$, then $d(f)=d(F)/N$, so we can just focus our attention on the fixed divisor of a primitive polynomial in $\Z[X]$.

We want to give another interpretation of the fixed divisor of a polynomial $f\in\Z[X]$ by considering the maximal ideals of $\IZ$ containing $f(X)$ and looking at their contraction to $\Z[X]$. We recall first the definition of unitary ideal given in \cite{McQ}.

\begin{Def}
 An ideal $I\subseteq\IZ$ is \textbf{unitary} if $I\cap\Z\not=0$.
\end{Def}
That is, an ideal $I$ of $\IZ$ is unitary if it contains a non-zero integer, or, equivalently, $I\Q[X]=\Q[X]$ (where $I\Q[X]$ denotes the extension ideal in $\Q[X]$). The whole ring $\IZ$ is clearly a principal unitary ideal generated by $1$.

The next results are probably well-known, but for the ease of the reader we report them. The first lemma says that a principal unitary ideal $I$ is generated by a non-zero integer, which generates the contraction of $I$ to $\Z$. In particular, this lemma establishes a bijective correspondence between the nonzero ideals of $\Z$ and the set of principal unitary ideals of $\IZ$. 
\vskip0.3cm
\begin{Lemma}\label{1}
Let $I\subseteq\IZ$ be a principal unitary ideal. If $I\cap\Z=n\Z$ with $n\not=0$ then $I=n\IZ$. In particular, $n\IZ\cap\Z=n\Z$. Moreover, $n_1\IZ=n_2\IZ$ with $n_1,n_2\in\Z$ if and only if $n_1=\pm n_2$.
\end{Lemma}
\Proof If $I=(f)$ for some $f\in\IZ$ then $\deg(f)=0$ since a non-zero integer $n$ is in $I$. Since $f(X)$ is integer-valued it must be equal to an integer and so it is contained in $I\cap\Z=n\Z$. Hence we get the first statement of the lemma. If $n_1\IZ=n_2\IZ$ then $n_1=n_2 f$ with $f\in\IZ$; this forces $f$ to be a non-zero integer, so that $n_1$ divides $n_2$. Similarly, we get that $n_2$ divides $n_1$. $\Box$
\vskip0.4cm
\begin{Lemma}\label{2}
Let $I_1,I_2\subseteq\IZ$ be principal unitary ideals. Then $I_1\cap I_2$ is a principal unitary ideal too.
\end{Lemma}
\Proof Suppose $I_i=n_i\IZ$, where $n_i\in\Z$, $n_i\Z=I_i\cap\Z$, for $i=1,2$. We have $n_1\Z\cap n_2\Z=n\Z$, where $n=\textnormal{lcm}\{n_1,n_2\}$. The ideal $I_1\cap I_2$ is unitary since $n\in I_1\cap I_2$. In particular, we have $I_1\cap I_2\supseteq n\IZ$. We have to prove that $I_1\cap I_2\subseteq n\IZ$. Let $f\in I_1\cap I_2$. Then $f(\Z)\subset n_1\Z\cap n_2\Z=n\Z$, so that $\frac{f(X)}{n}\in\IZ$. $\Box$
\vskip0.6cm
The previous lemma implies the following decomposition for a principal unitary ideal generated by an integer $n$, with prime factorization $n=\prod_i p_i^{a_i}$. We have
$$n\IZ=\bigcap_i p_i^{a_i}\IZ=\prod_i p_i^{a_i}\IZ$$
where the last equality holds because the ideals $p_i^{a_i}\Z$ are coprime in $\Z$, hence they are coprime in $\IZ$. 

We are now ready to give the following definition.
\begin{Def}\label{fixdivIZ}
Let $f\in\IZ$. The \textbf{extended fixed divisor} of $f(X)$ is the minimal ideal of the set $\{n\IZ\;|\;n\in\Z,f\in n\IZ\}$.
We denote this ideal by $D(f)$.

\end{Def}
\vskip0.3cm

Equivalently, in the above definition, we require that $n\IZ$ contains the principal ideal in $\IZ$ generated by the polynomial $f(X)$. Lemma \ref{1} and \ref{2} show that the minimal ideal in the above definition does exist: it is equal to the intersection of all the principal unitary ideals containing $f(X)$. Notice that the extended fixed divisor is an ideal of $\IZ$, while the fixed divisor is an ideal of $\Z$. The polynomial $f(X)$ is image primitive if and only if its extended fixed divisor is the whole ring $\IZ$. In the next sections we will study the extended fixed divisor by considering the $p$-part of it, namely the principal unitary ideals of the form $p^n\IZ$, $p\in\Z$ being prime and $n$ a positive integer. 

The following proposition gives a link between the fixed divisor and the extended fixed divisor: the latter is the extension of the former and conversely. So each of them gives information about the other one.
\vskip0.3cm
\begin{Prop}
Let $f\in\IZ$. Then we have:
\begin{itemize}
 \item[a)] $D(f)\cap\Z=d(f)$
 \item[b)] $d(f)\IZ=D(f)$
\end{itemize}
\end{Prop}
\Proof Let $d,D\in\Z$ be such that $d(f)=d\Z$ and $D(f)=D\IZ$. Since $d(f)\IZ=d\IZ$ is a principal unitary ideal containing $f(X)$, from the definition of extended fixed divisor, we have $D(f)\subseteq d\IZ$. In particular, $D\geq d$. We also have $f(X)/D\in\IZ$ and so $d\geq D$, by characterization 1) of Lemma \ref{fixdiv}. Hence we get a).
From that we deduce that $d(f)\subseteq D(f)$, so statement b) follows. $\Box$
\vskip0.5cm
As already remarked in \cite{CaCh}, the rings $\Z$ and $\IZ$ share the same units, namely $\{\pm1\}$. Then \cite[Proposition 2.1]{CaCh} can be restated as follows.
\vskip0.3cm
\begin{Prop}[Cahen-Chabert]
Let $f\in\IZ$ be irreducible in $\Q[X]$. Then $f(X)$ is irreducible in $\IZ$ if and only if $f(X)$ is not contained in any proper principal unitary ideal of $\IZ$.
\end{Prop}
\vskip0.1cm
The next lemma has been given in \cite{ChMcCl} and is analogous to the Gauss Lemma for polynomials in $\Z[X]$ which are irreducible in $\IZ$.
\begin{Lemma}[Chapman-McClain]
Let $f\in\Z[X]$ be a primitive polynomial. Then $f(X)$ is irreducible in $\IZ$ if and only if it is irreducible in $\Z[X]$ and image primitive.
\end{Lemma}
\vskip0.3cm

For example, the polynomial $f(X)=X^2+X+2$ is irreducible in $\Q[X]$ and also in $\Z[X]$ since it is primitive (because of Gauss Lemma). But it is reducible in $\IZ$ since its extended fixed divisor is not trivial, namely it is the ideal $2\IZ$. So in $\IZ$ we have the following factorization:
$$f(X)=2\cdot\frac{X^2+X+2}{2}$$
and indeed this is a factorization into irreducibles in $\IZ$, since the latter polynomial is image primitive and irreducible in $\Q[X]$, and by \cite[Lemma 1.1]{CaCh},  the irreducible elements in $\Z$ remain irreducible in $\IZ$. So the study of the extended fixed divisor of the elements in $\IZ$ is a first step toward studying the factorization of the elements in this ring (which is not a unique factorization domain).
\vskip0.5cm
Here is an overview of the content of the paper. At the beginning of the next section we recall the structure of the prime spectrum of $\IZ$. Then, for a fixed prime $p$, we describe the contractions to $\Z[X]$ of the maximal unitary ideals of $\IZ$ containing $p$ (Lemma \ref{Mp}). In Theorem \ref{pIZ} we describe the ideal $I_p$ of $\Z[X]$ of those polynomials whose fixed divisor is divisible by $p$, namely the contraction to $\Z[X]$ of the principal unitary ideal $p\IZ$, which is the ideal of integer-valued polynomials whose extended fixed divisor is contained in $p\IZ$. It turns out that $I_p$ is the intersection of the aforementioned contractions. In the third section we generalize the result of the second section to prime powers, by means of a structure theorem of Loper regarding unitary ideals of $\IZ$. We consider the contractions to $\Z[X]$ of the powers of the prime unitary ideals of $\IZ$ (Lemma \ref{Mpn}). In Remark \ref{primarydecpnRZ} we give a description of the structure of the set of these contractions; that allows us to give the primary decomposition of the ideal $I_{p^n}=p^n\IZ\cap\Z[X]$, made up of those polynomials whose fixed divisor is divisible by a prime power $p^n$. We shall see that we have to distinguish two cases: $p\leq n$ and $p>n$ (see also the examples in Remark \ref{Controesempi}). In Theorem \ref{pnIZ} we describe $I_{p^n}$ in the case $p\leq n$. This result was already known in a slightly different context by Dickson (see \cite[p. 22, Theorem 27]{Dickson}), but our different proof uses the primary decomposition of $I_{p^n}$ and that gives an insight to generalize the result to the second case. In Proposition \ref{Qp0p>n} we give a set of generators for the primary components of $I_{p^n}$, in the case $p>n$. Finally in the last section, as an application, we explicitly compute the ideal $I_{p^{p+1}}$.

\vskip2cm

\section{Fixed divisor via Spec($\IZ$)}

The study of the prime spectrum of the ring $\IZ$ began in \cite{Chab}. We recall that the prime ideals of $\IZ$ are divided into two different categories, unitary and non-unitary. Let $P$ be a prime ideal of $\IZ$. If it is unitary then its intersection with the ring of integers is a principal ideal generated by a prime $p$.\\
\vskip0.2cm
\noindent\textbf{Non-unitary prime ideals: } $P\cap\Z=\{0\}$.\\
\noindent In this case $P$ is a prime (non-maximal) ideal and it is of the form
$$\frak{B}_q=q\Q[X]\cap\IZ$$
for some $q\in\Q[X]$ irreducible. By Gauss Lemma we may suppose that $q\in\Z[X]$ is irreducible and primitive.
\vskip0.2cm
\noindent\textbf{Unitary prime ideals: } $P\cap\Z=p\Z$.\\
\noindent In this case $P$ is maximal and is of the form
$$\frak{M}_{p,\alpha}=\{f\in\IZ\;|\;f(\alpha)\in p\Z_p\}$$
for some $p$ prime in $\Z$ and some $\alpha\in\Z_p$, the ring of $p$-adic integers. We have $\frak{M}_{p,\alpha}=\frak{M}_{q,\beta}$ if and only if $(p,\alpha)=(q,\beta)$. So if we fix the prime $p$, the elements of $\Z_p$ are in bijection with the unitary prime ideals of $\IZ$ above the prime $p$. Moreover, $\frak{M}_{p,\alpha}$ is height $1$ if and only if $\alpha$ is transcendental over $\Q$. If $\alpha$ is algebraic over $\Q$ and $q(X)$ is its minimal polynomial then $\frak{M}_{p,\alpha}\supset\frak{B}_q$. We have $\frak{B}_q\subset\frak{M}_{p,\alpha}$ if and only if $q(\alpha)=0$. Every prime ideal of $\IZ$ is not finitely generated.

For a detailed study of Spec($\IZ$) see \cite{CaCh0}.
\vskip0.4cm
If we denote by $d(f,\Z_p)$ the fixed divisor of $f\in\IZ$ viewed as a polynomial over the ring of $p$-adic integers $\Z_p$ (that is, $d(f,\Z_p)$ is the ideal $(f(\alpha)\,|\,\alpha\in\Z_p)$), Gunji and McQuillan in \cite{GMcQ} observed that
\begin{equation*}
d(f)=\bigcap_{p}d(f,\Z_p)
\end{equation*}
where the intersection is taken over the set of primes in $\Z$. Moreover, $d(f,\Z_p)=d(f)\Z_p\subset\Z_p$. Remember that given an ideal $I\subset\Z$ and a prime $p$ we have $I\Z_p=\Z_p$ if and only if $I\not\subset (p)$, so that in the previous equation we have a finite intersection. Since $\Z_p$ is a DVR we have $d(f,\Z_p)=p^n\Z_p$, for some integer $n$ (which of course depends on $p$), so that the exact power of $p$ which divides $f(\Z)$ is the same as the power of $p$ dividing $f(\Z_p)$.
Without loss of generality, we can restrict our attention to the $p$-part of the fixed divisor of a polynomial $f\in\Z[X]$. We begin our research by finding those polynomials in $\Z[X]$ whose fixed divisor is divisible by a fixed prime $p$, namely the ideal $p\IZ\cap\Z[X]$.
\vskip0.5cm
\begin{Lemma}\label{Mp}
Let $p$ be a prime and $\alpha\in\Z_p$. Then $\frak{M}_{p,\alpha}\cap\Z[X]=(p,X-a)$, where $a\in\Z$ is such that $\alpha\equiv a\pmod p$. Moreover, if $\beta\in\Z_p$ is another $p$-adic integer, we have $\frak{M}_{p,\alpha}\cap\Z[X]=\frak{M}_{p,\beta}\cap\Z[X]$ if and only if $\alpha\equiv\beta\pmod p$.
\end{Lemma}
\Proof Let $a$ be an integer as in the statement of the lemma; it exists since $\Z$ is dense in $\Z_p$ for the $p$-adic topology. We immediately see that $p$ and $X-a$ are in $\frak{M}_{p,\alpha}$. Then the conclusion follows since $(p,X-a)$ is a maximal ideal of $\Z[X]$ and $\frak{M}_{p,\alpha}\cap\Z[X]$ is not equal to the whole ring $\Z[X]$. The second statement follows from the fact that $(p,X-a)=(p,X-b)$ if and only if $a\equiv b\pmod p$. $\Box$
\vskip0.6cm
We have just seen that the contraction of $\frak{M}_{p,\alpha}$ to $\Z[X]$ depends only on the residue class modulo $p$ of $\alpha$. So, if $p$ is a fixed prime, the contractions of $\frak{M}_{p,\alpha}$ to $\Z[X]$ as $\alpha$ ranges through $\Z_p$ are made up of $p$ distinct maximal ideals, namely 
$$\{\frak{M}_{p,\alpha}\cap\Z[X]\;|\;\alpha\in\Z_p\}=\{(p,X-j)\,|\,j\in\{0,\ldots,p-1\}\}.$$ 
Conversely, the set of prime ideals of $\IZ$ above a fixed maximal ideal of the form $(p,X-j)$ is $\{\,\frak{M}_{p,\alpha}\,|\,\alpha\in\Z_p\,,\,\alpha\equiv j\pmod p\}$, since $\frak{B}_q$ are non-unitary ideals and $p$ is the only prime integer in $\frak{M}_{p,\alpha}$. 
 
For a prime $p$ and an integer $j\in\{0,\ldots,p-1\}$, we set:
$$\mathcal{M}_{p,j}=\mathcal{M}_{j}\doteqdot(p,X-j).$$
Whenever the notation $\mathcal{M}_{p,j}$ is used, it will be implicit that $j\in\{0,\ldots,p-1\}$.

The next lemma computes the intersection of the ideals $\mathcal{M}_{p,j}$, for a fixed prime $p$, by finding an ideal whose primary decomposition is given by this intersection (and its primary components are precisely the $p$ ideals $\mathcal{M}_{p,j}$).  From now on we will omit the index $p$.
\begin{Lemma}\label{pRZx}
Let $p\in\Z$ be a prime.  Then we have
$$\bigcap_{j=0,\ldots,p-1}\mathcal{M}_{j}=\left(p,\prod_{j=0,\ldots,p-1}(X-j)\right).$$
\end{Lemma}
\Proof Let $J$ be the ideal on the right-hand side. If $P$ is a prime minimal over $J$, then we see immediately that $P=\mathcal{M}_{j}$ for some $j\in\{0,\ldots,p-1\}$, since $\mathcal{M}_{j}$ is a maximal ideal. Conversely, every such a maximal ideal contains $J$ and is minimal over it. Then the minimal primary decomposition of $J$ is of the form
$$J=\bigcap_{j=0,\ldots,p-1}Q_j$$
where $Q_j$ is an $\mathcal{M}_{j}$-primary ideal. Since $X-i\not\in\mathcal{M}_{j}$ for all $i\in\{0,\ldots,p-1\}\setminus\{j\}$, we have $(X-j)\in Q_j$, so indeed $Q_j=(p,X-j)$ for each $j=0,\ldots,p-1$. $\Box$

\vskip0.7cm
The next proposition characterizes the principal unitary ideals in $\IZ$ generated by a prime $p$.
\begin{Prop}\label{pR}
 Let $p\in\Z$ be a prime. Then the principal unitary ideal $p\IZ$ is equal to
 $$p\IZ=\bigcap_{\alpha\in\Z_p}\frak{M}_{p,\alpha}.$$
\end{Prop}
\Proof We trivially have that $p\IZ$ is contained in the above intersection, since $p$ is in every ideal of the form $\frak{M}_{p,\alpha}$. On the other hand, this intersection is equal to $\{f\in\IZ|f(\Z_p)\subset p\Z_p\}$. If $f(X)$ is in this intersection, since $f(X)$ is integer-valued and $p\Z_p\cap\Z=p\Z$, we have $f(\Z)\subset p\Z$. This is equivalent to saying that $f(X)/p\in\IZ$, that is, $f\in p\IZ$. $\Box$\\
\vskip0.2cm
In particular, the previous proposition implies that $\IZ$ does not have the finite character property (we recall that a ring has this property if every non-zero element is contained in a finite number of maximal ideals).

From the above results we get the following theorem, which characterizes the ideal of polynomials with integer coefficients whose fixed divisor is divisible by a prime $p$, that is, the ideal $p\IZ\cap\Z[X]$. 
\begin{Theorem}\label{pIZ}
Let $p$ be a prime. Then
$$p\IZ\cap\Z[X]=\left(p,\prod_{j=0,\ldots,p-1}(X-j)\right).$$
\end{Theorem}
\vskip0.4cm
Notice that Lemma \ref{pRZx} gives the primary decomposition of $p\IZ\cap\Z[X]$, so $\mathcal{M}_{j}$ for $j=0,\ldots,p-1$ are exactly the prime ideals belonging to it. As a consequence of this theorem we get the following well-known result: if $f\in\Z[X]$ is primitive and $p$ is a prime such that $d(f)\subseteq p$ then $p\leq\deg(f)$. This immediately follows from the theorem, since the degree of $\prod_{j=0,\ldots,p-1}(X-j)$ is $p$.

We remark that by Fermat's little theorem the ideal on the right-hand side of the statement of Theorem \ref{pIZ} is equal to $\left(p,X^p-X\right)$. This amounts to saying that the two polynomials $X\cdot\ldots\cdot(X-(p-1))$ and $X^p-X$ induce the same polynomial function on $\Z/p\Z$.

\vskip1cm

\section{Contraction of primary ideals}\label{contrprimary}
\vskip0.5cm
We remark that Proposition \ref{pR} also follows from a general result contained in \cite{Loper}: every unitary ideal in $\IZ$ is an intersection of powers of unitary prime ideals (namely the maximal ideals  $\frak{M}_{p,\alpha}$). In particular, every $\frak{M}_{p,\alpha}$-primary ideal is a power of $\frak{M}_{p,\alpha}$ itself, since $\frak{M}_{p,\alpha}$ is maximal. From the same result we also have the following characterization of the powers of $\frak{M}_{p,\alpha}$, for any positive integer $n$:
$$\frak{M}_{p,\alpha}^n=\{f\in\IZ\;|\;f(\alpha)\in p^n\Z_p\}.$$
This fact implies the following expression for the principal unitary ideal generated by $p^n$:
\begin{equation}\label{pnR}
p^n\IZ=\bigcap_{\alpha\in\Z_p}\frak{M}_{p,\alpha}^n.
\end{equation}
We remark again that the previous ideal is made up of those integer-valued polynomials whose extended fixed divisor is contained in $p^n\IZ$. Similarly to the previous case $n=1$ (see Theorem \ref{pIZ}) we want to find the contraction of this ideal to $\Z[X]$, in order to find the polynomials in $\Z[X]$ whose fixed divisor is divisible by $p^n$. We set:
\begin{equation}\label{Ipn}
I_{p^n}\doteqdot p^n\IZ\cap\Z[X].
\end{equation}
Notice that by (\ref{pnR}) we have $I_{p^n}=\bigcap_{\alpha\in\Z_p}(\frak{M}_{p,\alpha}^n\cap\Z[X])$.

Like before, we begin by finding the contraction to $\Z[X]$ of $\frak{M}_{p,\alpha}^n$, for each $\alpha\in\Z_p$. The next lemma is a generalization of Lemma \ref{Mp}.
\vskip0.4cm
\begin{Lemma}\label{Mpn}
Let $p$ be a prime, $n$ a positive integer and $\alpha\in\Z_p$. Then $\frak{M}_{p,\alpha}^n\cap\Z[X]=(p^n,X-a)$, where $a\in\Z$ is such that $\alpha\equiv a\pmod{p^n}$. The ideal $\frak{M}_{p,\alpha}^n\cap\Z[X]$ is  $\mathcal{M}_{p,j}$-primary, where $j\equiv\alpha\pmod p$. Moreover, if $\beta\in\Z_p$ is another $p$-adic integer, we have $\frak{M}_{p,\alpha}^n\cap\Z[X]=\frak{M}_{p,\beta}^n\cap\Z[X]$ if and only if $\alpha\equiv\beta\pmod{p^n}$.
\end{Lemma}
\Proof The case $n=1$ has been done in Lemma \ref{Mp}. For the general case, let $a\in\Z$ be such that $a\equiv\alpha\pmod{p^n}$ (again, such an integer exists since $\Z$ is dense in $\Z_p$ for the $p$-adic topology). We have $(p^n,X-a)\subset\frak{M}_{p,\alpha}^n\cap\Z[X]$ (notice that if $n>1$ then $(p^n,X-a)$ is not a prime ideal). To prove the other inclusion let $f\in\frak{M}_{p,\alpha}^n\cap\Z[X]$. By the Euclidean algorithm in $\Z[X]$ (the leading coefficient of $X-a$ is a unit) we have 
\begin{equation*}
f(X)=q(X)(X-a)+f(a)
\end{equation*}
Since $f(\alpha)\in p^n\Z_p$ and $p^n|a-\alpha$ we have $p^n|f(a)$. Hence, $f\in(p^n,X-a)$ as we wanted. Since $\frak{M}_{p,\alpha}^n$ is an $\frak{M}_{p,\alpha}$-primary ideal in $\IZ$ and the contraction of a primary ideal is a primary ideal, by Lemma \ref{Mp} we get the second statement. Finally, like in the proof of Lemma \ref{Mp}, we immediately see that $(p^n,X-a)=(p^n,X-b)$ if and only if $a\equiv b\pmod{p^n}$, which gives the last statement of the lemma. $\Box$ \\
\vskip0.6cm
\begin{Oss}\label{divresto}
It is worth to write down the fact that we used in the above proof: given a polynomial $f\in\Z[X]$, we have
\begin{equation}\label{**}
f\in(p^n,X-a)\Longleftrightarrow f(a)\equiv0\pmod{p^n}
\end{equation}
\end{Oss}
\vskip0.9cm

\begin{Oss}\label{primarydecpnRZ}
If $p$ is a fixed prime and $n$ is a positive integer, the Lemma \ref{Mpn} implies
$$\mathcal{I}_{p,n}\doteqdot\{\frak{M}_{p,\alpha}^n\cap\Z[X]\;|\;\alpha\in\Z_p\}=\{(p^n,X-i)\;|\;i=0,\ldots,p^n-1\}.$$
Let us consider an ideal $I=\frak{M}_{p,\alpha}^n\cap\Z[X]=(p^n,X-i)$ in $\mathcal{I}_{p,n}$, with $i\in\Z$, $i\equiv\alpha\pmod{p^n}$. It is quite easy to see that $I$ contains $(\frak{M}_{p,\alpha}\cap\Z[X])^n=\mathcal{M}_{p,j}^n=(p,X-j)^n$, where $j\in\{0,\ldots,p-1\}$, $j\equiv\alpha\pmod{p}$ (notice that $j\equiv i\pmod p$). If $n>1$ this containment is strict, since $X-i\not\in(p,X-j)^n$. 
We can group the ideals of $\mathcal{I}_{p,n}$ according to their radical: there are $p$ radicals of these $p^n$ ideals, namely the maximal ideals $\mathcal{M}_{p,j}$, $j=0,\ldots,p-1$. This amounts to making a partition of the residue classes modulo $p^n$ into $p$ different sets of elements congruent to $j$ modulo $p$, for $j=0,\ldots,p-1$; each of these sets has cardinality $p^{n-1}$. Correspondingly we have:
$$\mathcal{I}_{p,n}=\bigcup_{j=0,\ldots,p-1}\mathcal{I}_{p,n,j}$$
where $\mathcal{I}_{p,n,j}\doteqdot\{(p^n,X-i)\;|\;i=0,\ldots,p^n-1,i\equiv j\pmod p\}$, for $j=0,\ldots,p-1$. Every ideal in $\mathcal{I}_{p,n,j}$ is $\mathcal{M}_{p,j}$-primary and it contains the $n$-th power of its radical, namely $\mathcal{M}_{p,j}^n$.

Now we want to compute the intersection of the ideals in $\mathcal{I}_{p,n}$, which is equal to the ideal $I_{p^n}$ in $\Z[X]$ (see (\ref{pnR}) and (\ref{Ipn})). We can express this intersection as an intersection of $\mathcal{M}_{p,j}$-primary ideals as we have said above, in the following way (in the first equality we make use of equation (\ref{pnR}) and Lemma \ref{Mpn}):
\begin{equation}\label{intIpn}
I_{p^n}=\bigcap_{i=0,\ldots,p^n-1}(p^n,X-i)=\bigcap_{j=0,\ldots,p-1}\mathcal{Q}_{p,n,j}
\end{equation}  
where 
$$\mathcal{Q}_{p,n,j}\doteqdot\bigcap_{i\equiv j(\textnormal{mod } p)}(p^n,X-i)$$
(notice that the intersection is taken over the set $\{i\in\{0,\ldots,p^n-1\} \mid i\equiv j\pmod p\}$). The ideal $\mathcal{Q}_{p,n,j}$ is an $\mathcal{M}_{p,j}$-primary ideal, for $j=0,\ldots,p-1$, since the intersection of $M$-primary ideals is an $M$-primary ideal. We will omit the index $p$ in $\mathcal{Q}_{p,n,j}$ and in $\mathcal{M}_{p,j}$ if that will be clear from the context. The $\mathcal{M}_{p,j}$-primary ideal $\mathcal{Q}_{n,j}$ is just the intersection of the ideals in $\mathcal{I}_{p,n,j}$, according to the partition we made. It is equal to the set of polynomials in $\Z[X]$ which modulo $p^n$ are zero at the residue classes congruent to $j$ modulo $p$ (see (\ref{**}) of Remark \ref{divresto}). We remark that (\ref{intIpn}) is the minimal primary decomposition of $I_{p^n}$. Notice that there are no embedded components in this primary decomposition, since the prime ideals belonging to it (the minimal primes containing $I_{p^n}$) are $\{\mathcal{M}_{j}\,|\,j=0,\ldots,p-1\}$, which are maximal ideals.
 
We recall that if $I$ and $J$ are two coprime ideals in a ring $R$, that is  $I+J=R$, then $IJ=I\cap J$ (in general only the inclusion $IJ\subset I\cap J$ holds). The condition for two ideals $I$ and $J$ to be coprime amounts to saying that $I$ and $J$ are not contained in a same maximal ideal $M$, that is, $I+J$ is not contained in any maximal ideal $M$. If $M_1$ and $M_2$ are two distinct maximal ideals then they are coprime, and the same holds for any of their respective powers. If $R$ is Noetherian, then every primary ideal $Q$ contains a power of its radical and moreover if the radical of $Q$ is maximal then also the converse holds (see \cite{North}). So if $Q_i$ is an $M_i$-primary ideal for $i=1,2$ and $M_1,M_2$ are distinct maximal ideals, then $Q_1$ and $Q_2$ are coprime.
 
Since $\{\mathcal{M}_{j}\}_{j=0,\ldots,p-1}$ are $p$ distinct maximal ideals, for what we have just said above we have 
\begin{equation*}
\bigcap_{j=0,\ldots,p-1}\mathcal{Q}_{n,j}=\prod_{j=0,\ldots,p-1}\mathcal{Q}_{n,j}.
\end{equation*}
\end{Oss}
\vskip0.3cm
Now we want to describe the $\mathcal{M}_{j}$-primary ideals $\mathcal{Q}_{n,j}$, for $j=0,\ldots,p-1$. The next lemma gives a relation of containment between these ideals and the $n$-th powers of their radicals.
\vskip0.3cm
\begin{Lemma}\label{Qpj} Let $p$ be a fixed prime and $n$ a positive integer. For each $j=0,\ldots,p-1$, we have 
$$\mathcal{Q}_{n,j}\supseteq\mathcal{M}_{j}^n.$$
\end{Lemma}
\Proof The statement follows from Remark \ref{primarydecpnRZ}. $\Box$
\vskip0.7cm 
As a consequence of this lemma, we get the following result:
\begin{Cor}\label{pnIZcontpprodn}
 Let $p$ be a fixed prime and $n$ a positive integer. Then we have:
 $$I_{p^n}\supseteq\Big(p,\prod_{j=0,\ldots,p-1}(X-j)\Big)^n.$$
\end{Cor}
\Proof By (\ref{intIpn}) and Lemma \ref{Qpj} we have
$$I_{p^n}=\prod_{j=0,\ldots,p-1}\mathcal{Q}_{n,j}\supseteq\prod_{j=0,\ldots,p-1}\mathcal{M}_{j}^n$$ 
where the last containment follows from Lemma \ref{Qpj}. Finally, by Lemma \ref{pRZx}, the product of the ideals $\mathcal{M}_{j}^n$ is equal to
$$\prod_{j=0,\ldots,p-1}\mathcal{M}_{j}^n=\Big(p,\prod_{j=0,\ldots,p-1}(X-j)\Big)^n$$
Notice that the product of the $\mathcal{M}_{j}$'s is actually equal to their intersection, since they are maximal coprime ideals. $\Box$
\vskip0.3cm
The last formula of the previous proof gives the primary decomposition of the ideal \newline $(p,\prod_{j=0,\ldots,p-1}(X-j))^n$.

\vskip0.7cm
\begin{Oss}\label{Controesempi}
\noindent In general, for a fixed $j\in\{0,\ldots,p-1\}$, the reverse containment of Lemma \ref{Qpj} does not hold, that is, the $n$-th power of $\mathcal{M}_{j}$ can be strictly contained in the $\mathcal{M}_{j}$-primary ideal $\mathcal{Q}_{n,j}$. For example (again, we use (\ref{**}) to prove the containment):
$$X(X-2)\in\left(\bigcap_{k=0,\ldots,3}(2^3,X-2k)\right)\setminus(2,X)^3$$
Because of that, in general we do not have an  equality in Corollary \ref{pnIZcontpprodn}. For example, let $p=2$ and $n=3$. We have
$$X(X-1)(X-2)(X-3)\in I_{2^3}\setminus(2,X(X-1))^3.$$
\noindent It is also false that 
$$\bigcap_{i=0,\ldots,p^n-1}(p^n,X-i)=\left(p^n,\prod_{i=0,\ldots,p^n-1}(X-i)\right).$$
See for example: $p=2$, $n=2$:  $2X(X-1)\in \bigcap_{i=0,\ldots,3}(4,X-i)\setminus\left(4,\prod_{i=0,\ldots,3}(X-i)\right)$.
\end{Oss}
\vskip0.4cm

We want to study under which conditions the ideal $\mathcal{Q}_{n,j}$ is equal to $\mathcal{M}_{j}^n$. Our aim is to find a set of generators for $\mathcal{Q}_{n,j}$. For  $f\in\mathcal{Q}_{n,j}$, we have $f\in(p^n,X-i)$ for each $i\equiv j \pmod p$, $i\in\{0,\ldots,p^n-1\}$. By (\ref{**}) that means $p^n|f(i)$ for each such an $i$. Equivalently, such a polynomial has the property that modulo $p^n$ it is zero at the $p^{n-1}$ residue classes of $\Z/p^n\Z$ which are congruent to $j$ modulo $p$. 

Without loss of generality, we proceed by considering the case $j=0$. We set $\mathcal{M}=\mathcal{M}_{0}=(p,X)$ and $\mathcal{Q}_{n}=\mathcal{Q}_{n,0}=\bigcap_{i\equiv0\pmod p}(p^n,X-i)$. Let $f\in \mathcal{Q}_{n}$, of degree $m$. We have
\begin{equation}\label{k=1}
f(X)=q_1(X)X+f(0)
\end{equation}
where $q_1\in\Z[X]$ has degree equal to $m-1$. Since $f\in(p^n,X)$ we have $p^n|f(0)$.

Since $f\in(p^n,X-p)$, we have $p^n|f(p)=q_1(p)p+f(0)$, so $p^{n-1}|q_1(p)$. By the Euclidean algorithm, 
\begin{equation}\label{k=1.1}
q_1(X)=q_2(X)(X-p)+q_1(p)
\end{equation}
for some polynomial $q_2\in\Z[X]$ of degree $m-2$. So
$$f(X)=q_2(X)(X-p)X+q_1(p)X+f(0).$$
We set $R_1(X)=q_1(p)X+f(0)$. Then $R_1\in\mathcal{M}^n$, since $p^{n-1}|q_1(p)$ and $p^n|f(0)$. Since $f\in(p^n,X-2p)$, we have $p^n|f(2p)=q_2(2p)2p^2+q_1(p)2p+f(0)$. If $p>2$ then $p^{n-2}|q_2(2p)$, because $p^{n}|q_1(p)2p+f(0)$. If $p=2$ then we can just say $p^{n-3}|q_2(2p)$. By the Euclidean algorithm again, we have
$$q_2(X)=q_3(X)(X-2p)+q_2(2p)$$
for some $q_3\in\Z[X]$. So we have
$$f(X)=q_3(X)(X-2p)(X-p)X+q_2(2p)(X-p)X+q_1(p)X+f(0).$$
Like before, if we set $R_2(X)=q_2(2p)(X-p)X+q_1(p)X+f(0)$, we have $R_2\in\mathcal{M}^n$ if $p>2$, or $R_2\in\mathcal{Q}_{n}$ if $p=2$.

We define now the following family of polynomials:
\begin{Def}\label{Gk}
 For each $k\in\N$, $k\geq1$, we set
 \begin{equation*}
G_{p,0,k}(X)=G_{k}(X)\doteqdot\prod_{h=0,\ldots,k-1}(X-hp).
\end{equation*}
We also set $G_0(X)\doteqdot1$.
\end{Def}

From now on, we will omit the index $p$ in the above notation. 

Notice that the polynomials $G_k(X)$, whose degree for each $k$ is $k$, enjoy these properties: 
\begin{itemize}\label{propGk}
\item[i)] For every $t\in\Z$, $G_k(tp)=p^k t(t-1)\ldots(t-(k-1))$. Hence, the highest power of $p$ which divides all the integers in the set $\{G_k(tp)\,|\,t\in\Z\}$ is $p^{k+v_p(k!)}$. It is easy to see that $k+v_p(k!)=v_p((pk)!)$.
\item[ii)] $G_k(X)=(X-kp)G_{k-1}(X).$
\item[iii)] since for every integer $h$, $X-hp\in\mathcal{M}$, we have $G_k(X)\in\mathcal{M}^k$. We remark that $k$ is the maximal integer with this property, since $\deg(G_k)=k$ and $G_k(X)$ is primitive (since monic). 
\end{itemize}

Recall that, by Lemma \ref{Qpj}, for every integer $n$ we have $\mathcal{Q}_n\supseteq\mathcal{M}^n$. By property iii) above $G_k\in\mathcal{M}^n$ if and only if $n\leq k$. By property i) we have $G_k\in\mathcal{Q}_n$ if and only if $k+v_p(k!)\geq n$.
From these remarks, it is very easy to deduce that, in the case $p\geq n$, if $G_k\in\mathcal{Q}_n$ then $G_k\in\mathcal{M}^n$. In fact, if that is not the case, it follows from above that $k<n$. Since $n\leq p$ we get $k+v_p(k!)=k$. Since $G_k\in\mathcal{Q}_n$, we have $n\leq k$, contradiction. 

The next lemma gives a sort of division algorithm between an element of $\mathcal{Q}_n$ and the polynomials $\{G_k(X)\}_{k\in\N}$. In particular, we will deduce that $\mathcal{Q}_n=\mathcal{M}^n$, if $p\geq n$.

\begin{Lemma}\label{qGR0+1}
 Let $p$ be a prime and $n$ a positive integer. Let $f\in\mathcal{Q}_{p,n,0}=\mathcal{Q}_{n}$ be of degree $m$. Then for each $1\leq k\leq m$ there exists $q_k\in\Z[X]$ of degree $m-k$ such that 
$$f(X)=q_k(X)G_k(X)+R_{k-1}(X)$$ 
where $R_{k-1}(X)\doteqdot\sum_{h=1,\ldots,k-1}q_h(hp)G_h(X)$ for $k\geq2$ and $R_0(X)\doteqdot f(0)$. We also have $q_k(X)=q_{k+1}(X)(X-kp)+q_k(kp)$ for $k=1,\ldots,m-1$.
Moreover, for each such a $k$ the following hold:
 \begin{itemize}
  \item[i)] $p^{n-v_p((pk)!)}|q_k(kp)$, if $v_p((pk)!)<n$.
  \item[ii)] $q_k(kp)G_k(X)\in\mathcal{Q}_n$ and if $k<p$ then $q_k(kp)G_k(X)\in\mathcal{M}^n$.
  \item[iii)] If $m\leq p$ then $R_{k-1}\in\mathcal{M}^n$ for $k=1,\ldots,m$.\\ If $m>p$ then $R_{k-1}\in\mathcal{M}^n$ for $k=1,\ldots,p$ and $R_{k-1}\in\mathcal{Q}_{n}$ for $k=p+1,\ldots,m$. 
 \end{itemize}
\end{Lemma}
\vskip0.5cm
\Proof We proceed by induction on $k$. The case $k=1$ follows from (\ref{k=1}), and by (\ref{k=1.1}) we have the last statement regarding the relation between $q_1(X)$ and $q_2(X)$. Suppose now the statement is true for $k-1$, so that 
$$f(X)=q_{k-1}(X)G_{k-1}(X)+R_{k-2}(X)$$
with $R_{k-2}(X)\doteqdot\sum_{h=1,\ldots,k-2}q_h(hp)G_h(X)$ and
\begin{itemize}
 \item[-] $p^{n-v_p((p(k-1))!)}|q_{k-1}((k-1)p)$, if $v_p((p(k-1)!))<n$,
 \item[-] $q_{k-1}((k-1)p)G_{k-1}(X)$ belongs to $\mathcal{Q}_n$ and if $k-1<p$ it belongs to $\mathcal{M}^n$,
 \item[-] $R_{k-2}\in\mathcal{Q}_{n}$ and if $k-2<p$ then $R_{k-2}\in\mathcal{M}^n$.
\end{itemize}
We divide $q_{k-1}(X)$ by $(X-(k-1)p)$ and we get 
$$q_{k-1}(X)=q_{k}(X)(X-(k-1)p)+q_{k-1}((k-1)p)$$
for some polynomial $q_k\in\Z[X]$ of degree $m-k$. We substitute this expression of $q_{k-1}(X)$ in the equation of $f(X)$ at the step $k-1$ and we get:
\begin{equation}\label{fk}
f(X)=q_{k}(X)(X-(k-1)p)G_{k-1}(X)+R_{k-1}(X),
\end{equation}
where $R_{k-1}(X)\doteqdot q_{k-1}((k-1)p)G_{k-1}(X)+R_{k-2}(X)$. This is the expression of $f(X)$ at step $k$, since $(X-(k-1)p)G_{k-1}(X)$ is equal to $G_k(X)$. By the inductive assumption, $R_{k-1}\in\mathcal{Q}_n$ and if $k-1<p$ we also have $R_{k-1}\in\mathcal{M}^n$. We still have to verify i) and ii).

We evaluate the expression (\ref{fk}) in $X=kp$ and we get
$$f(kp)=q_k(kp)G_k(kp)+R_{k-1}(kp)=q_k(kp)p^k k!+R_{k-1}(kp).$$
Since $p^n$ divides both $f(kp)$ and $R_{k-1}(kp)$ (by definition of $\mathcal{Q}_n$), if $v_p((pk)!)<n$ we get that $q_k(kp)$ is divisible by $p^{n-v_p((pk)!)}$,  which is statement i) at the step $k$. Notice that $q_k(kp)G_k(X)$ is zero modulo $p^n$ on every integer congruent to zero modulo $p$; hence, $q_k(kp)G_k(X)\in\mathcal{Q}_n$. Moreover, $k<p\Leftrightarrow v_p(k!)=0$, so in that case  $q_k(kp)G_k(X)\in\mathcal{M}^n$. So ii) follows. $\Box$

\vskip0.8cm

Notice that by formula (\ref{**}) of Remark \ref{divresto}, under the assumptions of Lemma \ref{qGR0+1}  we have for each $k\in\{1,\ldots,p-1\}$ that
$$q_k\in(p^{n-k},X-kp)$$
(see i) of Lemma \ref{qGR0+1}: in this case $v_p((pk)!)=k$). If $k=m=\deg(f)$ then $q_k\in\Z$. Hence, we get the following expression for a polynomial $f\in\mathcal{Q}_{n}$ in the case $p\geq n>m$ (this assumption is not restrictive, since $X^n\in\mathcal{Q}_n$):
\begin{equation}\label{***}
f(X)=q_m G_m(X)+R_{m-1}(X)=q_m G_m(X)+\sum_{k=1,\ldots, m-1}q_k(kp)G_k(X)
\end{equation}
where $q_m\in\Z$ is divisible by $p^{n-m}$ and $R_{m-1}(X)$ is in $\mathcal{M}^n$.

The next proposition determines the primary components $\mathcal{Q}_{n,j}$ of $I_{p^n}$ of (\ref{intIpn}) in the case $p\geq n$. It shows that in this case the containment of Lemma \ref{Qpj} is indeed an equality.
\vskip0.6cm
\begin{Prop}\label{Qp0p<=n}
Let $p\in\Z$ be a prime and $n$ a positive integer such that $p\geq n$. Then for each $j=0,\ldots,p-1$ we have 
$$\mathcal{Q}_{n,j}=\mathcal{M}_{j}^n.$$
\end{Prop}
\Proof It is sufficient to prove the statement for $j=0$: for the other cases we consider the $\Z[X]$-automorphisms $\pi_j(X)=X-j$, for $j=1,\ldots,p-1$, which permute the ideals $\mathcal{Q}_{n,j}$ and $\mathcal{M}_{j}$. Let $\mathcal{Q}_{n}=\mathcal{Q}_{n,0}$ and $\mathcal{M}=\mathcal{M}_{0}$. 

The inclusion $(\supseteq)$ follows from Lemma \ref{Qpj}. For the other inclusion $(\subseteq)$, let $f(X)$ be in $\mathcal{Q}_{n}$. We can assume that the degree $m$ of $f(X)$ is less than $n$, since $X^n$ is the smallest monic monomial in $\mathcal{Q}_{n}$. By equation (\ref{***}) above, $f(X)$ is in $\mathcal{M}^n$, since $p^{n-m}$ divides $q_m$, $G_m\in\mathcal{M}^m$ and $R_{m-1}\in \mathcal{M}^n$ by Lemma \ref{qGR0+1} (notice that $m-1<p$). $\Box$
\vskip0.5cm
\begin{Oss} In the case $p\geq n$, Lemma \ref{qGR0+1} implies that $\mathcal{Q}_{n}$ is generated by $\{p^{n-m}G_m(X)\}_{0\leq m\leq n}$: it is easy to verify that these polynomials are in $\mathcal{Q}_{n}$ (using (\ref{**}) again) and (\ref{***}) implies that every polynomial $f\in\mathcal{Q}_n$ is a $\Z$-linear combination of $\{p^{n-m}G_m(X)\}_{0\leq m\leq n}$, since $q_{m}(mp)$ is divisible by $p^{n-m}$, for each of the relevant $m$.
\end{Oss}

The following theorem gives a description of the ideal $I_{p^n}$ in the case $p\geq n$. In this case the containment of the Corollary \ref{pnIZcontpprodn} becomes an equality.

\vskip0.5cm
\begin{Theorem}\label{pnIZ}
Let $p\in\Z$ be a prime and $n$ a positive integer such that $p\geq n$. Then the ideal in $\Z[X]$ of those polynomials whose fixed divisor is divisible by $p^n$ is equal to
\begin{equation*}
I_{p^n}=\Big(p,\prod_{i=0,\ldots,p-1}(X-i)\Big)^n.
\end{equation*}
\end{Theorem}
\Proof 
By Proposition \ref{Qp0p<=n}, for each $j=0,\ldots,p-1$ the ideal $\mathcal{Q}_{n,j}$ is equal to $\mathcal{M}_j^n$.
So, by the last formula of the proof of Corollary \ref{pnIZcontpprodn}, we get the statement. $\Box$
\vskip0.8cm
As a consequence, we have the following remark. Let $p$ be a prime and $n$ a positive integer less than or equal to $p$. Let $f\in I_{p^n}$ such that the content of $f(X)$ is not divisible by $p$. Then $\deg(f)\geq np$, since $np=\deg(\prod_{i=0,\ldots,p-1}(X-i)^n)$. Another well-known result in this context is the following: if we fix the degree $d$ of such a polynomial $f$, then the maximum $n$ such that $f\in I_{p^n}$ is bounded by $n\leq\sum_{k\geq1} [d/p^k] =v_p(d!)$. 

If we drop the assumption $p\geq n$, the ideal $\mathcal{Q}_{n,j}$ may strictly contain $\mathcal{M}_{j}^n$, as we observed  in Remark \ref{Controesempi}. The next proposition shows that this is always the case, if $p<n$. This result follows from Lemma \ref{qGR0+1} as Proposition \ref{Qp0p<=n} does, and it covers the remaining case $p<n$. It is stated for the case $j=0$. Remember that $\mathcal{M}=(p,X)$ and $\mathcal{Q}_{n}=\bigcap_{i\equiv0\pmod p}(p^n,X-i)$.
\begin{Prop}\label{Qp0p>n}
 Let $p\in\Z$ be a prime and $n$ a positive integer such that $p<n$. Then we have 
$$\mathcal{Q}_{n}=\mathcal{M}^n+(q_{n,p}G_p(X),\ldots,q_{n,n-1}G_{n-1}(X))$$
where, for each $k=p,\ldots,n-1$, $q_{n,k}$ is an integer defined as follows:
\begin{equation*}
q_{n,k} \doteqdot \left\{
\begin{array}{rl}
p^{n-v_p((pk)!)} & \text{, if } v_p((pk)!)<n\\
1 & \text{, otherwise}
\end{array} \right.
\end{equation*}
In particular, $\mathcal{M}^n$ is strictly contained in $\mathcal{Q}_{n}$.
\end{Prop}
\Proof We begin by proving the containment $(\supseteq)$. Lemma \ref{Qpj} gives $\mathcal{M}^n\subseteq\mathcal{Q}_{n}$. We have to show that the polynomials $q_{n,k}G_k(X)$, for $k\in\{p,\ldots,n-1\}$, lie in $\mathcal{Q}_{n}$. This follows from property i) of the polynomials $G_k(X)$ and the definition of $q_{n,m}$.

Now we prove the other containment $(\subseteq)$. Let $f\in\mathcal{Q}_{n}$ be of degree $m$. If $m<p$ then $f\in\mathcal{M}^n$ (see Lemma \ref{qGR0+1} and in particular (\ref{***})). So we suppose $p\leq m$. By Lemma \ref{qGR0+1} we have
\begin{equation}\label{finQ}
f(X)=\sum_{k=p,\ldots,m}q_k(kp)G_k(X)+R_{p-1}(X)
\end{equation}
where $R_{p-1}(X)=\sum_{k=1,\ldots,p-1}q_k(kp)G_k(X)\in\mathcal{M}^n$ and $q_{m}\in\Z$, so that $q_m(mp)=q_{n,m}$. Then, since $q_{n,k}=p^{n-v_p((pk)!)}|q_k(kp)$ if $v_p((pk)!)<n$, it follows that the first sum on the right-hand side of the previous equation belongs to the ideal $(q_{n,p}G_p(X),\ldots,q_{n,n-1}G_{n-1}(X))$. 
For the last sentence of the proposition, we remark that the polynomials $\{q_{n,k}G_k(X)\}_{k=p,\ldots,n-1}$ are not contained in $\mathcal{M}^n$: in fact, for each $k\in\{p,\ldots,n-1\}$, by property iii) of the polynomials $G_k(X)$ we have that the minimal integer $N$ such that $q_{n,k}G_k(X)$ is contained in $\mathcal{M}^N$  is $n-v_p(k!)$   if $v_p((pk)!)=k+v_p(k!)<n$ and it is $k$ otherwise. In both cases it is strictly less than $n$ (since $v_p(k!)\geq1$, if $k\geq p$). $\Box$
\vskip0.5cm
\begin{Oss}
The following remark allows us to obtain another set of generators for  $\mathcal{Q}_{n}$. We set
\begin{equation}\label{m}
\overline{m}=\overline{m}(n,p)\doteqdot\min\{m\in\N\;|\;v_p((pm)!)\geq n\}
\end{equation}
Remember that $v_p((pm)!)=m+v_p(m!)$. If $p\geq n$ then $\overline{m}=n$ and if $p<n$ then $p\leq \overline{m}< n$.

Suppose now $p<n$. Then for each $m\in\{\overline{m},\ldots,n\}$ we have $v_p((pm)!)\geq n$, since the function $e(m)=m+v_p(m!)$ is increasing. So for each such $m$ we have $q_{n,m}=1$, hence $G_m\in(G_{\overline{m}}(X))$. So we have the equalities:
\begin{align}\label{Qn}
\mathcal{Q}_{n}&=\mathcal{M}^n+(q_{n,m}G_m(X)\;|\;m=p,\ldots,\overline{m})\\
               &=(q_{n,m}G_m(X)\;|\;m=0,\ldots,\overline{m})\nonumber
\end{align}
where $q_{n,m}=p^{n-m}$, for $m=0,\ldots,p-1$, and for $m=p,\ldots,\overline{m}$ is defined as in the statement of Proposition \ref{Qp0p>n}. The containment $(\supseteq)$ is just an easy verification using the properties of the polynomials $G_m(X)$; the other containment follows by (\ref{finQ}). 
\end{Oss}
\vskip0.5cm
We can now group together Proposition \ref{Qp0p<=n} and \ref{Qp0p>n} into the following one:
\begin{Prop}\label{Qp0}
 Let $p\in\Z$ be a prime and $n$ a positive integer. Then we have 
$$\mathcal{Q}_{n}=(q_{n,0}G_0(X),\ldots,q_{n,\overline{m}}G_{\overline{m}}(X))$$
where $\overline{m}=\min\{m\in\N\;|\;v_p((pm)!)\geq n\}$ and for each $m=0,\ldots,\overline{m}$, $q_{n,m}$ is an integer defined as follows:
\begin{equation*}
q_{n,m} \doteqdot \left\{
\begin{array}{rl}
p^{n-v_p((pm)!)} &, m<\overline{m}\\
1 &, m=\overline{m}
\end{array} \right.
\end{equation*}
\end{Prop}
\vskip0.5cm
It is clear what the primary ideals $\mathcal{Q}_j$, for $j=1,\ldots,p-1$, look like:
\begin{align*}
\mathcal{Q}_{n,j}=\bigcap_{i\equiv j(\textnormal{mod } p)}(p^n,X-i)&=\mathcal{M}_{j}^n+(q_{n,p}G_p(X-j),\ldots,q_{n,\overline{m}}G_{\overline{m}}(X-j))\\
    &=(q_{n,0} G_0(X-j),\ldots,q_{n,\overline{m}}G_{\overline{m}}(X-j))
\end{align*}
In fact, for each $j=1,\ldots,p-1$, it is sufficient to consider the automorphisms of $\Z[X]$ given by $\pi_j(X)=X-j$. It is straightforward to check that $\pi_j(I_{p^n})=I_{p^n}$. Moreover, $\pi(\mathcal{Q}_{n,0})=\mathcal{Q}_{n,j}$ and $\pi(\mathcal{M}_{0})=\mathcal{M}_{j}$ for each such a $j$, so that $\pi_j$ permutes the primary components of the ideal $I_{p^n}$.

The ideal $I_{p^n}=p^n\IZ\cap\Z[X]$ was studied in \cite{Band} in a slightly different context, as the kernel of the natural map 
$\varphi_n:\Z[X]\to\Phi_n$, where the latter is the set of functions from $\Z/p^n\Z$ to itself. In that article a recursive formula is given for a set of generators of this ideal. Our approach gives a new point of view to describe this ideal. 

For other works about the ideal $I_{p^n}$ in a slightly different context, see \cite{KellOls}, \cite{Lewis}, \cite{NivWar}. This ideal is important in the study of the problem of the polynomial representation of a function from $\Z/p^n\Z$ to itself.

\vskip1cm
\section{Case $I_{p^{p+1}}$}
As a corollary we give an explicit expression for the ideal $I_{p^n}$ in the case $n=p+1$. By Proposition \ref{Qp0p>n} the primary components of $I_{p^{p+1}}$ are
\begin{equation}\label{Qp+1j}
\mathcal{Q}_{p+1,j}=\mathcal{M}_j^{p+1}+(G_{p}(X-j))
\end{equation}
for $j=0,\ldots,p-1$.
\begin{Cor}
 $$I_{p^{p+1}}=\Big(p,\prod_{i=0,\ldots,p-1}(X-i)\Big)^{p+1}+(H(X))$$
 where $H(X)=\prod_{i=0,\ldots,p^2-1}(X-i)$.
\end{Cor}
\vskip0.3cm
We want to stress that the polynomial $H(X)$ is not contained in the first ideal of the right-hand side of the statement. In \cite{Band} a similar result is stated with another polynomial $H_2(X)$  instead of our $H(X)$. Indeed the two polynomials, as already remarked in \cite{Band}, are congruent modulo the ideal $(p,\prod_{i=0,\ldots,p-1}(X-i))^{p+1}$.
\vskip0.3cm
\Proof Like before, we set $\mathcal{Q}_{p,p+1,j}=\mathcal{Q}_{p+1,j}$. The containment $(\supseteq)$ follows from corollary \ref{pnIZcontpprodn} and because the polynomial $H(X)$ is equal to $\prod_{j=0,\ldots,p-1}G_p(X-j)$ and for each $j=0,\ldots,p-1$ the polynomial $G_p(X-j)$ is in $\mathcal{Q}_{p+1,j}$ by Proposition \ref{Qp0p>n}. Since $\mathcal{Q}_{p+1,j}$, for $j=0,\ldots,p-1$, are exactly the primary components of $I_{p^{p+1}}$ (see (\ref{intIpn})), we get the claim.

Now we prove the other containment $(\subseteq)$. Let $f\in I_{p^{p+1}}=\bigcap_{j=0,\ldots,p-1}\mathcal{Q}_{p+1,j}$. By (\ref{Qp+1j}) we have: 
$$f(X)\equiv C_{p,j}(X)G_p(X-j)\pmod{\mathcal{M}_{j}^{p+1}}$$
for some $C_{p,j}\in\Z[X]$, for $j=0,\ldots,p-1$.

Since the ideals $\{\mathcal{M}_{j}^{p+1}=(p,X-j)^{p+1}\,|\,j=0,\ldots,p-1\}$ are pairwise coprime (because they are powers of distinct maximal ideals, respectively), by the Chinese Remainder Theorem we have the following isomorphism:
\begin{equation}\label{CRT}
\Z[X]/\left(\prod_{j=0}^{p-1} \mathcal{M}_{j}^{p+1}\right)\cong\Z[X]/\mathcal{M}_0^{p+1}\times\ldots\times\Z[X]/\mathcal{M}_{p-1}^{p+1}
\end{equation}

We need now the following lemma, which tells us what is the residue of the polynomial $H(X)$ modulo each ideal $\mathcal{M}_j^{p+1}$:
\vskip0.5cm
\begin{Lemma}\label{resHp+1}
 Let $p$ be a prime and let $H(X)=\prod_{j=0,\ldots,p-1}G_p(X-j)$. Then for each $k=0,\ldots,p-1$ we have
 $$H(X)\equiv-G_p(X-k)\pmod{\mathcal{M}_{k}^{p+1}}$$
\end{Lemma}
\Proof Let $k\in\{0,\ldots,p-1\}$ and set $I_k=\{0,\ldots,p-1\}\setminus\{k\}$. For each $j\in I_k$ we have 
$G_p(k-j)\equiv(k-j)^p\pmod p$.
We have
$$H(X)+G_p(X-k)=G_p(X-k)[1+\prod_{j\in I_k}G_p(X-j)]$$
Since $G_p(X-k)\in\mathcal{M}_{k}^{p}$ we have just to prove that $T_k(X)=1+\prod_{j\in I_k}G_p(X-j)\in\mathcal{M}_{k}$. By formula (\ref{**}) in Remark \ref{divresto} it is sufficient to prove that $T_k(k)$ is divisible by $p$. We have
\begin{align*}
T_k(k)&\equiv1+\prod_{j\in I_k}(k-j)^p\pmod p\\
    &\equiv1+(\prod_{s=1,\dots , p-1}s)^p\pmod p\\
    &\equiv1+(p-1)!^p\pmod p\\
    &\equiv(1+(p-1)!)^p\pmod p
\end{align*}
which is congruent to zero by Wilson's theorem. $\Box$
\vskip0.3cm We finish now the proof of the corollary.

By the Chinese Remainder Theorem, there exists a polynomial $P\in\Z[X]$ such that $P(X)\equiv-C_{p,j}(X)\pmod{\mathcal{M}_j^{p+1}}$, for each $j=0,\ldots,p-1$. Then by the previous lemma $P(X)H(X)\equiv C_{p,j}(X)G_p(X-j)\pmod{\mathcal{M}_j^{p+1}}$ and so again by the isomorphism (\ref{CRT}) above we have
$$f(X)\equiv P(X)H(X)\pmod{\prod_{j=0,\ldots,p-1} \mathcal{M}_j^{p+1}}$$
so we are done since $\prod_{j=0,\ldots,p-1} \mathcal{M}_j^{p+1}=(p,\prod_{i=0,\ldots,p-1}(X-i))^{p+1}$ (see the proof of Corollary \ref{pnIZcontpprodn}). $\Box$

\vskip1.1cm
\subsection*{\textbf{Acknowledgments}}
I want to thank the referee for his/her suggestions. I wish also to thank Sophie Frisch for the useful references. I deeply thank Elvira Carlotti for her continuous and strenuous help.  The author was supported by the Austrian Science Foundation (FWF), Project Number P23245-N18.


\vskip0.5cm
\addcontentsline{toc}{section}{Bibliography}

\begin{thebibliography}{99}
\bibliographystyle{plain}
\bibitem{ACaChS}  D. F. Anderson, P.-J. Cahen, S.T. Chapman and W. Smith, \emph{Some factorization properties of the ring of integer-valued polynomials.}, Zero-dimensional commutative rings (Knoxville, TN, 1994), 125-142, Lecture Notes in Pure and Appl. Math., 171, Dekker, New York, 1995.
\bibitem{Band} A. Bandini. \emph{Functions $f:\Z/p^n\Z\to\Z/p^n\Z$ induced by polynomials of $\Z[X]$}. Ann. Mat. Pura Appl. (4) 181 (2002), no.1, pp. 95-104.
\bibitem{Chab} J.-L. Chabert, \emph{Anneaux de polyn{\^o}mes {\`a} valeurs enti{\`e}res et
anneaux de Fatou}, Bull. Soc. Math. France 99 (1971), 273-283.
\bibitem{CaCh0} P.-J. Cahen, J.-L. Chabert, \emph{Integer-valued polynomials}, Math. Surv. Monogr., vol. 48, Amer. Math. Soc., Providence, 1997.
\bibitem{CaCh} P.-J. Cahen, J.-L. Chabert, \emph{Elasticity for integral-valued polynomials.}, J. Pure Appl. Algebra 103 (1995), no. 3, 303-311. 
\bibitem{ChMcCl}  S.T.Chapman, B. McClain. \emph{Irreducible polynomials and full elasticity in rings of integer-valued polynomials}, J. Algebra 293 (2005), no. 2, 595-610.
\bibitem{Dickson} L. E. Dickson. \emph{Introduction to the theory of numbers.} Univ. Chicago Press, Chicago, 1929.
\bibitem{GMcQ} H. Gunji and D. L. McQuillan. \emph{On a class of ideals in an algebraic number field}. J. Number Theory 2 (1970) 207--222.
\bibitem{KellOls} G. Keller and F.R. Olson. \emph{Counting polynomial functions (mod $p^n$)}. Duke Math. J. 35 (1968), 835-838.
\bibitem{Lewis} D. J. Lewis. \emph{Ideals and Polynomial Functions},  Amer. J. Math. 78 (1956), 71-77.
\bibitem{Loper}  K. A. Loper. \emph{Ideals of integer-valued polynomial rings}, Comm. Algebra 25 (1997), no. 3, 833-845.
\bibitem{McQ} D. L. McQuillan. \emph{On Pr\"ufer domains of polynomials}. J. reine angew. Math. 358 (1985), 162-178.
\bibitem{NivWar} I. Niven, L. J. Warren. \emph{A generalization of Fermat's theorem}. Proc. Amer. Math. Soc. 8 (1957), 306-313.
\bibitem{North} D. G. Northcott. \emph{Ideal Theory}. Cambridge University Press, 1953.

\end{thebibliography}
\end{document}